\documentclass[11pt, twoside]{article}
\usepackage{latexsym}
\usepackage{amsmath}
\usepackage{amssymb}
\usepackage[all]{xy}
\usepackage{amsfonts}
\usepackage{verbatim}
\usepackage{amsthm}
\usepackage{mathrsfs}
\usepackage{epsfig}
\usepackage{xy}
\usepackage{tensor}
\usepackage{array}
\usepackage{stmaryrd}
\usepackage{graphicx,color}
\usepackage{xcolor}
\usepackage[colorlinks=true,linkcolor=blue,citecolor=blue]{hyperref}
\usepackage{tikz}
\usetikzlibrary{arrows,calc}
\usepackage{etex}
\usepackage{mathdots}
\usepackage{float}
\usepackage{graphics}
\usepackage{pdflscape}
\usepackage{extarrows}
\usepackage{anysize,hyperref}
\input xypic
\xyoption{all}
\usepackage{bm}
\usepackage[perpage,symbol]{footmisc}
\usepackage{setspace}
\SelectTips{cm}{10}

\renewcommand{\cal}{\mathcal}

\def\C{\mathscr{C}}

\def\E{\mathbb{E}}
\def\s{\mathfrak{s}}

\def\del{\delta}
\def\dr{\ar@{->}[r]}

\def\Hom{\mbox{Hom}}

\begin{document}
\baselineskip=15pt
\title{\Large{\bf Abelian hearts of twin cotorsion pairs on extriangulated categories\footnotetext{ Panyue Zhou was supported by the National Natural Science Foundation of China (Grant No. 11901190) and the Scientific Research Fund of Hunan Provincial Education Department (Grant No. 19B239).} }}
\medskip
\author{Qiong Huang and Panyue Zhou}

\date{}

\maketitle
\def\blue{\color{blue}}
\def\red{\color{red}}

\newtheorem{theorem}{Theorem}[section]
\newtheorem{lemma}[theorem]{Lemma}
\newtheorem{corollary}[theorem]{Corollary}
\newtheorem{proposition}[theorem]{Proposition}
\newtheorem{conjecture}{Conjecture}
\theoremstyle{definition}
\newtheorem{definition}[theorem]{Definition}
\newtheorem{question}[theorem]{Question}
\newtheorem{notation}[theorem]{Notation}
\newtheorem{remark}[theorem]{Remark}
\newtheorem{remark*}[]{Remark}
\newtheorem{example}[theorem]{Example}
\newtheorem{example*}[]{Example}

\newtheorem{construction}[theorem]{Construction}
\newtheorem{construction*}[]{Construction}

\newtheorem{assumption}[theorem]{Assumption}
\newtheorem{assumption*}[]{Assumption}

\baselineskip=17pt
\parindent=0.5cm

\begin{abstract}
\baselineskip=16pt
It was shown recently that the heart of a twin cotorsion pair
on an extriangulated category is semi-abelian. In this article,
we consider a special kind of hearts of twin cotorsion pairs induced by $d$-cluster tilting subcategories in extriangulated categories.
We give a necessary and sufficient
condition for such hearts to be abelian. In particular,
  we also can see that such hearts are hereditary. As an application, this generalizes the work by Liu in an exact case, thereby providing new insights in a triangulated case. \\[0.5cm]
\medskip
\end{abstract}
\pagestyle{myheadings}
\markboth{\rightline {\scriptsize Q. Huang and P. Zhou\hspace{2mm}}}
         {\leftline{\scriptsize  Abelian hearts of twin cotorsion pairs on extriangulated categories}}

\section{Introduction}
The notion of cotorsion pairs was first introduced by Scale in \cite{S}, it has been studied in the representation theory.
 A careful look reveals that what is necessary to define a cotorsion pair on a category is the existence of an ${\rm Ext}^1$ bifunctor with appropriate properties. Nakaoka and Palu \cite{NP}
formalized the notion of an \emph{extriangulated category} by extracting those properties of ${\rm Ext}^1$ on exact categories and on triangulated categories that seem relevant from the point-of-view of cotorsion pairs.
Exact categories
and triangulated categories are examples of extriangulated categories, while there are some other examples of
extriangulated categories which are neither exact nor triangulated, see \cite{NP,ZZ,HZZ}. Hence,
many results on exact categories and triangulated categories can be unified in the same framework.
Nakaoka and Palu also defined (twin) cotorsion
pairs on extriangulated categories in \cite{NP}.
They showed that the heart of a twin cotorsion pair was always semi-abelian and the heart of a single cotorsion pair became abelian.
This unified the similar arguments on (twin) cotorsion pairs in exact categories \cite{L1} and triangulated categories \cite{N}.

Recently, Liu \cite{L2} studied a special kind of hearts of twin cotorsion pairs induced by $d$-cluster tilting subcategories in exact categories.
He gave a sufficient-necessary condition when such hearts become abelian.
We know that extriangulated categories are a simultaneous generalization of exact categories and
triangulated categories. A very natural question is whether these similar results hold in extriangulated categories.
In this article, we give an affirmative answer of this question.

This article is organized as follows: In Section 2, we review some elementary definitions
that we need to use, extriangulated categories, twin cotorsion pairs and $d$-cluster tilting subcategories. In Section 3, we prove our main result, see Theorem \ref{b6}.

\section{Preliminaries}
We briefly recall some definitions and basic properties of extriangulated categories from \cite{NP}.
We omit some details here, but the reader can find them in \cite{NP}.

Let $\mathcal{C}$ be an additive category equipped with an additive bifunctor
$$\mathbb{E}: \mathcal{C}^{\rm op}\times \mathcal{C}\rightarrow {\rm Ab},$$
where ${\rm Ab}$ is the category of abelian groups. For any objects $A, C\in\mathcal{C}$, an element $\delta\in \mathbb{E}(C,A)$ is called an $\mathbb{E}$-extension.
Let $\mathfrak{s}$ be a correspondence which associates an equivalence class $$\mathfrak{s}(\delta)=\xymatrix@C=0.8cm{[A\ar[r]^x
 &B\ar[r]^y&C]}$$ to any $\mathbb{E}$-extension $\delta\in\mathbb{E}(C, A)$. This $\mathfrak{s}$ is called a {\it realization} of $\mathbb{E}$, if it makes the diagrams in \cite[Definition 2.9]{NP} commutative.
 A triplet $(\mathcal{C}, \mathbb{E}, \mathfrak{s})$ is called an {\it extriangulated category} if it satisfies the following conditions.
\begin{itemize}
\item $\mathbb{E}\colon\mathcal{C}^{\rm op}\times \mathcal{C}\rightarrow \rm{Ab}$ is an additive bifunctor.

\item $\mathfrak{s}$ is an additive realization of $\mathbb{E}$.

\item $\mathbb{E}$ and $\mathfrak{s}$  satisfy the compatibility conditions in \cite[Definition 2.12]{NP}.
\end{itemize}

We collect the following terminology from \cite{NP}.

\begin{definition}
Let $(\C,\E,\s)$ be an extriangulated category.
\begin{itemize}
\item[(1)] A sequence $A\xrightarrow{~x~}B\xrightarrow{~y~}C$ is called a {\it conflation} if it realizes some $\E$-extension $\del\in\E(C,A)$.
    In this case, $x$ is called an {\it inflation} and $y$ is called a {\it deflation}.

\item[(2)] If a conflation  $A\xrightarrow{~x~}B\xrightarrow{~y~}C$ realizes $\delta\in\mathbb{E}(C,A)$, we call the pair $( A\xrightarrow{~x~}B\xrightarrow{~y~}C,\delta)$ an {\it $\E$-triangle}, and write it in the following way.
$$A\overset{x}{\longrightarrow}B\overset{y}{\longrightarrow}C\overset{\delta}{\dashrightarrow}$$
We usually do not write this $``\delta"$ if it is not used in the argument.

\item[(3)] Let $A\overset{x}{\longrightarrow}B\overset{y}{\longrightarrow}C\overset{\delta}{\dashrightarrow}$ and $A^{\prime}\overset{x^{\prime}}{\longrightarrow}B^{\prime}\overset{y^{\prime}}{\longrightarrow}C^{\prime}\overset{\delta^{\prime}}{\dashrightarrow}$ be any pair of $\E$-triangles. If a triplet $(a,b,c)$ realizes $(a,c)\colon\delta\to\delta^{\prime}$, then we write it as
$$\xymatrix{
A \ar[r]^x \ar[d]^a & B\ar[r]^y \ar[d]^{b} & C\ar@{-->}[r]^{\del}\ar[d]^c&\\
A'\ar[r]^{x'} & B' \ar[r]^{y'} & C'\ar@{-->}[r]^{\del'} &}$$
and call $(a,b,c)$ a {\it morphism of $\E$-triangles}.

\item[(4)] An object $P\in\C$ is called {\it projective} if
for any $\E$-triangle $\xymatrix{A\ar[r]^{x}&B\ar[r]^{y}&C\ar@{-->}[r]^{\delta}&}$ and any morphism $c\in\C(P,C)$, there exists $b\in\C(P,B)$ satisfying $yb=c$.
We denote the full subcategory of projective objects in $\C$ by $\cal P$. Dually, the full subcategory of injective objects in $\C$ is denoted by $\cal I$.

\item[(5)] We say $\C$ {\it has enough projectives}, if
for any object $C\in\C$, there exists an $\E$-triangle
$$\xymatrix{A\ar[r]^{x}&P\ar[r]^{y}&C\ar@{-->}[r]^{\delta}&}$$
satisfying $P\in\cal P$.  We can define the notion of having enough injectives dually.
\end{itemize}
\end{definition}

{\bf From now on to the end of the article}, let $k$ be a field and $(\C,\E,\s)$
be a Krull-Schmidt, Hom-finite, $k$-linear extriangulated category
with enough projectives $\mathcal{P}$ and enough injectives $\mathcal{I}$. When we say that $\mathcal A$ is a subcategory of a category, we always mean that $\mathcal A$ is full and is closed under isomorphisms, direct sums and direct summands.

The following notions can be found in \cite{NP} and \cite{LN}.
\begin{definition}
Let $\mathcal{C}'$ and $\mathcal{C}''$ be two subcategories of $\C$.

\begin{itemize}
\item[{\rm (a)}] Denote by CoCone( $\left.\mathcal{C}', \mathcal{C}''\right)$ the subcategory

$\left\{A \in \C \mid\right.$ there exists an $\mathbb{E}$-triangle $\xymatrix{A\ar[r]&C'\ar[r]&C''\ar@{-->}[r]&}$, $C' \in \mathcal{C}'$ and $\left.C'' \in \mathcal{C}''\right\}$.

\item[{\rm (b)}] Denote by Cone $\left(\mathcal{C}', \mathcal{C}''\right)$ the subcategory

$\left\{A \in \C \mid\right.$ there exists an $\mathbb{E}$-triangle $\xymatrix{C'\ar[r]&C''\ar[r]&A\ar@{-->}[r]&}$, $C' \in \mathcal{C}'$ and $\left.C'' \in \mathcal{C}''\right\}$.

\item[{\rm (c)}] Let $\Omega^{0} \mathcal{C}'=\mathcal{C}'$ and $\Omega \mathcal{C}'=\operatorname{CoCone}\left(\mathcal{P}, \mathcal{C}'\right),$ then we can define $\Omega^{i} \mathcal{C}'$ inductively:
$$
\Omega^{i} \mathcal{C}^{\prime}=\operatorname{CoCone}\left(\mathcal{P}, \Omega^{i-1} \mathcal{C}^{\prime}\right).
$$
We write an object $B$ in the form $\Omega C$ if it admits an $\mathbb{E}$-triangle $\xymatrix{B\ar[r]&P\ar[r]&C\ar@{-->}[r]&}$ where $P \in \mathcal{P}$.

\item[{\rm (d)}] Let $\Sigma^{0} \mathcal{C}'=\mathcal{C}', \Sigma \mathcal{C}'=\operatorname{Cone}\left(\mathcal{C}', \mathcal{I}\right),$ then we can define $\Sigma^{i} \mathcal{C}'$ inductively:
$$
\Sigma^{i} \mathcal{C}'=\operatorname{Cone}\left(\Sigma^{i-1} \mathcal{C}', \mathcal{I}\right).
$$
We write an object $B'$ in the form $\Sigma C'$ if it admits an $\mathbb{E}$-triangle $\xymatrix{C'\ar[r]&I\ar[r]&B'\ar@{-->}[r]&}$ where $I \in \mathcal{I}$.
\end{itemize}
\end{definition}

\begin{definition}{\cite[Definition 4.1 and Definition 4.12]{NP}}
Let $\mathcal{U}$ and $\mathcal{V}$ be two subcategories of $\C$. We call $(\mathcal{U}, \mathcal{V})$ a \emph{cotorsion pair} if it satisfies the following conditions.
\begin{itemize}
\item[{\rm (a)}] $\mathbb{E}(\mathcal{U}, \mathcal{V})=0$.
\item[{\rm (b)}] For any object $C \in \mathscr{C},$ there are two $\E$-triangles $$\xymatrix{V_{C}\ar[r]&U_{C}\ar[r]&C\ar@{-->}[r]&}~\mbox{and}~ \xymatrix{C\ar[r]&V^{C}\ar[r]&U^{C}\ar@{-->}[r]&}$$
satisfying $U_{C}, U^{C}\in U$ and $V_{C}, V^{C}\in V$.
\end{itemize}

A cotorsion pair $(\mathcal{U}, \mathcal{V})$ is said to be \emph{rigid} if it satisfies $\E(\mathcal U,\mathcal U)=0$.

Let $(\mathcal{S}, \mathcal{T})$ and $(\mathcal{U}, \mathcal{V})$ be cotorsion pairs on $\C$. Then the pair $(\mathcal{S}, \mathcal{T}), (\mathcal{U}, \mathcal{V})$ is called a \emph{twin cotorsion pair} if it satisfies $\mathbb{E}(\mathcal{S}, \mathcal{V})=0$, or equivalently $\mathcal S\subseteq\mathcal U$.
\end{definition}

\begin{remark}
Let $(\mathcal{U}, \mathcal{V})$ be a cotorsion pair on $\C$. Then
\begin{itemize}
\item[{\rm (a)}] $C\in \mathcal{U}$ if and only if $\mathbb{E}(C, \mathcal{V})=0$.
\item[{\rm (b)}] $C\in \mathcal{V}$ if and only if $\mathbb{E}(\mathcal{U}, C)=0$.
\item[{\rm (c)}] $\mathcal{U}$ and $\mathcal{V}$ are extension closed.
\item[{\rm (d)}] $\mathcal{P}\subseteq \mathcal{U}$ and $\mathcal{I}\subseteq \mathcal{V}$.
\end{itemize}
\end{remark}

\begin{definition}{\cite[Definition 2.5 and Definition 2.6]{LN}}
For any twin cotorsion pair $(\mathcal{S}, \mathcal{T}),(\mathcal{U}, \mathcal{V}),$ put $\mathcal{W}=\mathcal{T} \cap \mathcal{U}$ and call it the \emph{core} of $(\mathcal{U}, \mathcal{V}).$ We give the following definition.
\begin{itemize}
\item[{\rm (a)}] $\C^{+}=\operatorname{Cone}(\mathcal{V}, \mathcal{W})$. Namely, $\C^{+}$ is defined to be the full subcategory of $\C,$ consisting of objects $C$ which admits an $\E$-triangle
$$
\xymatrix{V_{C}\ar[r]&W_{C}\ar[r]&C\ar@{-->}[r]&}
$$
where $W_{C} \in \mathcal{W}$ and $V_{C} \in \mathcal{V}$. Clearly, we have $\mathcal{V}\subseteq \mathcal{T}\subseteq \C^{+}$.
\item[{\rm (b)}] $\C^{-}=\mathrm{CoCone}(\mathcal{W}, \mathcal{S})$. Namely, $\C^{-}$ is defined to be the full subcategory of $\C,$ consisting of objects $C$ which admits an $\E$-triangle
$$
\xymatrix{C\ar[r]&W^{C}\ar[r]&S^{C}\ar@{-->}[r]&}
$$
where $W^{C} \in \mathcal{W}$ and $S^{C} \in \mathcal{S}$. Clearly, we have $\mathcal{S}\subseteq \mathcal{U}\subseteq \C^{-}$.
\item[{\rm (c)}] $\mathcal{H}/ \mathcal{W}=(\C^{+} \cap \C^{-})/ \mathcal{W}$. The additive quotient $\mathcal{H}/ \mathcal{W}$ is called the \emph{heart} of $(\mathcal{S}, \mathcal{T}),(\mathcal{U}, \mathcal{V}).$
\end{itemize}
\end{definition}

Liu and Nakaoka \cite[Proposition 5.2]{LN} defined higher extension groups in an extriangulated category with enough projectives and enough injectives as $$\mathbb{E}^{i+1}(X, Y):= \mathbb{E}(\Omega^{i}X, Y)\cong \mathbb{E}(X, \Sigma^{i}Y)$$ for any $i \geq 0$.

\begin{definition}{\cite[Definition 5.3]{LN}}
Let $\mathcal{M}$ be a subcategory of $\C$.
\begin{itemize}
\item $\mathcal{M}$ is called \emph{d-cluster tilting} if it satisfies the following conditions:
\begin{itemize}
\item[{\rm (1)}] $\mathcal{M}$ is functorially finite in $\C$,
\item[{\rm (2)}] $X \in \mathcal{M}$ if and only if $\mathbb{E}^{i}(X, \mathcal{M})=0$ for any $i\in \{1, 2, \cdot \cdot \cdot , d-1\}$,
\item[{\rm (3)}] $X \in \mathcal{M}$ if and only if $\mathbb{E}^{i}(\mathcal{M}, X)=0$ for any $i\in \{1, 2, \cdot \cdot \cdot , d-1\}$.
\end{itemize}
\end{itemize}
\begin{itemize}
\item $\mathcal{M}$ is called \emph{n-rigid} if $\mathbb{E}^{i}(\mathcal{M}, \mathcal{M})=0$ for any $i\in \{1, 2, \cdot \cdot \cdot , n\}$.
\end{itemize}
\end{definition}

\section{Abelian hearts of twin cotorsion pairs}
\setcounter{equation}{0}
In this section, we extend the results in \cite{L2} to extriangulated categories.
In order to prove our main result, we need some preparations as follows.

\begin{lemma}\label{b1}
Let $A$ be an indecomposable object in $\C$.
\begin{itemize}
\item[{\rm (a)}] In the $\mathbb{E}$-triangle $\xymatrix{\Omega A\ar[r]^{g}&P\ar[r]^{f}&A\ar@{-->}[r]&}$ where $f$ is a right minimal, if $\mathbb{E}(A, P)=0$, then $\Omega A$ is indecomposable.
\item[{\rm (b)}] In the $\mathbb{E}$-triangle $\xymatrix{A\ar[r]^{k}&I\ar[r]^{h}&\Sigma A\ar@{-->}[r]&}$ where $k$ is a left minimal, if $\mathbb{E}(I, A)=0$, then $\Sigma A$ is indecomposable.
\end{itemize}
\end{lemma}

\proof The proof is an adaption of \cite[Lemma 2.1]{L2}.
We only prove (a), the proof of (b) is similar.

Assume $\Omega A=B_{1}\oplus B_{2}$ where $0\neq B_{1}$ is indecomposable. We have an $\mathbb{E}$-triangle of the form $$
B_1\oplus B_2\xrightarrow{~\left[\begin{smallmatrix}g_{1}&g_{2}\end {smallmatrix}\right]~}P\xrightarrow{~f~}A\xymatrix{\ar@{-->}[r]^{\delta}&}.$$
By \cite[Remark 2.16]{NP}, we know that inflations are closed under composition. Then $g_{1}$ is an inflation since $\left[\begin{smallmatrix}g_{1}&g_{2}\end {smallmatrix}\right]$ and $\left[\begin{smallmatrix}1\\0\end {smallmatrix}\right]$ are inflations. Thus it admits an $\mathbb{E}$-triangle $$\xymatrix{B_1\ar[r]^{g_{1}}&P\ar[r]^{f_1}&A_1\ar@{-->}[r]^{\sigma}&}.$$
Applying the functor $\Hom_{\C}(-,P)$ to the above $\E$-triangle, we have the following exact sequence:
$${\C}(P,P)\xrightarrow{~{\C}(g_1,~P)~}{\C}(B_1\oplus B_2,P)\xrightarrow{}\E(A,P)=0.$$
So there exists a morphism $p\colon P\to P$ such that
$p\left[\begin{smallmatrix}g_{1}&g_{2}\end {smallmatrix}\right]=g_1\left[\begin{smallmatrix}1&0\end {smallmatrix}\right]$.
Hence we obtain a morphism of $\mathbb{E}$-triangles
$$\xymatrix@R=1cm{B_{1}\oplus B_{2} \ar[r]^{\quad \left[\begin{smallmatrix}g_{1}&g_{2}\end {smallmatrix}\right]} \ar[d]^{\left[\begin{smallmatrix} 1& 0\end {smallmatrix}\right]} &P\ar[r]^{f} \ar[d]^{ p} & A\ar@{-->}[r]^{\delta} \ar[d]^{y}&\\
B_{1}\ar[r]^{g_1}\ar[d]^{\left[\begin{smallmatrix}1\\[1mm]0\end {smallmatrix}\right]} &P \ar[r]^{f_1}\ar@{=}[d] & A_1\ar[d]^x \ar@{-->}[r]^{\sigma}&\\
B_{1}\oplus B_{2} \ar[r]^{\quad \left[\begin{smallmatrix}g_{1}&g_{2}\end {smallmatrix}\right]}&P\ar[r]^{f} & A\ar@{-->}[r]^{\delta} &}$$
Since $A$ is indecomposable and $\C$ is Krull-Schmidt, $\textrm{End}_{\C}(A)$ is a local ring. Then, $xy$ is nilpotent or is an isomorphism.

If  $xy$ is nilpotent, i.e. there exists a positive integer $n$ such that $(xy)^{n}=0$.
We consider the following commutative diagram
$$\xymatrix{B_{1}\ar[r]^{g_{1}} \ar@{=}[d] &P\ar[r]^{f_{1}}\ar[d]^{p^{n+1}} & A_1\ar@{-->}[r]^{\sigma} \ar[d]^{x(xy)^{n}y=0}&\\
B_{1}\ar[r]^{g_{1}}&P \ar[r]^{f_1} & A_1 \ar@{-->}[r]^{\sigma}&}$$
 By \cite[Corollary 3.5]{NP}, we have that $1_{B_{1}}$ factors through $g_{1}$. Then $B_{1}$ is a direct summand of $P$. We write $P= B_{1}\oplus Q$. We have an $\mathbb{E}$-triangle
$$B_{1}\oplus B_{2}\xrightarrow{~\left[\begin{smallmatrix}g_{11}&g_{12} \\ g_{21}&g_{22}\end{smallmatrix}\right]~}B_{1}\oplus Q\xrightarrow{~ \left[\begin{smallmatrix}m&n\end {smallmatrix}\right]~}A\xymatrix{\ar@{-->}[r]^{\delta}&}$$
with $g_{11}$ is an isomorphism and $\left[\begin{smallmatrix}m&n\end {smallmatrix}\right]= f$. We have an isomorphism $B_{1}\oplus Q\xrightarrow{~\left[\begin{smallmatrix}g_{11}&0 \\ g_{21}&1_Q\end {smallmatrix}\right]~}B_{1}\oplus Q$ such that $\left[\begin{smallmatrix}m&n\end {smallmatrix}\right]\left[\begin{smallmatrix}g_{11}&0 \\ g_{21}&1_Q\end {smallmatrix}\right]= \left[\begin{smallmatrix}0&n\end {smallmatrix}\right]$, which is a contradiction since $f$ is right minimal.
Hence $xy$ is an isomorphism.
Note that $fp=xyf$, since $f$ is right minimal, we have that
$p$ is an isomorphism. It follows that
$B_1\oplus B_2\simeq B_1$, as desired.  \qed
\medskip

The notion of fully rigid subcategories in an extriangulated category was introduced by Liu and Zhou \cite{LZ2}, it can be regarded as a simultaneous generalization of \cite[Definition 5.1]{B} and \cite[Definition 1.5]{L2}.
We recall it here.

\begin{definition}{\cite[Definition 2.6]{LZ2}}\label{fullyrigid}
A subcategory $\mathcal{M}$ of $\C$ is called \emph{fully rigid} if
\begin{itemize}
\item[{\rm (a)}] it admits a rigid cotorsion pair $(\mathcal{M}, \mathcal{Y})$ and $\mathcal{M}\neq \mathcal{P}$,
\item[{\rm (b)}] any indecomposable object in $\C$ either belongs to $\mathcal{Y}$ or belongs to $\mathcal{H}:= \textrm{CoCone}(\mathcal{M}, \mathcal{M})$.

\end{itemize}
\end{definition}

\begin{remark}
From \cite[Remark 1.5]{LZ2}, we can see that when $\mathcal{M}$ is fully rigid, it is contravariantly finite, rigid, and $\C/  {\mathcal{M}}^{\perp _{1}}$ is an abelian category where ${\mathcal{M}}^{\perp _{1}}=\{ X\in \C \mid \mathbb{E}(\mathcal{M}, X)=0\}.$
 \end{remark}

 \begin{remark}
 From the definition of full rigid subcategories, we know that cluster tilting subcategories are fully rigid, but fully rigid subcategory are not necessarily cluster tilting, see \cite[Example 5.1]{LZ2}. Furthermore, $d$-cluster tilting is not always fully rigid when $d\geq 3$, see \cite[Example 3.1]{L2}.
 \end{remark}

Liu and Nakaoka showed that a $d$-cluster tilting subcategory $\mathcal{M}$ always admits a twin cotorsion pair $(\mathcal{M},\mathcal{Y}),(\mathcal{Y},\mathcal{N})$, see \cite[Theorem 5.14]{LN}. The heart of such a twin cotorsion pair is just $\C/ \mathcal{Y}$. We denote $\C/ \mathcal{Y}$ by $\underline{\C}$, for any morphism $f\in \C(X, Y)$, we denote its image in $\underline{\C}(X, Y)$ by $\underline{f}$. From now on to the end of the article,  let $(\mathcal{M}, \mathcal{Y})$ be a rigid  cotorsion pair. Then the heart of $(\mathcal{M},\mathcal{Y})$ is $\mathcal{H}/ \mathcal{M}$. Moreover, for a fully rigid subcategory $\mathcal{M}$, $$\underline{\C}\simeq \mathcal{H}/ \mathcal{M}\simeq \textrm{mod} (\Omega \mathcal{M}/ \mathcal{P})\simeq \textrm{mod} (\mathcal{M}/ \mathcal{P}),$$ and $\underline{\C}$ has enough projectives $\Omega \mathcal{M}/ \mathcal{P}$, see \cite[Proposition 2.11]{LZ2} and \cite[Theorem 4.10]{LN}.

Let $\mathcal{W}$ and $\mathcal{R}$ be subcategories of $\C$. We denote
$$
\mathcal{W}_{\mathcal{R}}=\{W\in \mathcal{W} \mid \textrm{W has no non-zero direct summands in}~\mathcal{R}\}\cup\{\mbox{zero objects of}~ \mathcal W\}.
$$

\begin{lemma}\label{h0}
If $(\mathcal{M},\mathcal{Y})$ is a rigid cotorsion pair, then $\mathcal{H}$ is closed under direct summands.
\end{lemma}

\proof This can be shown in the same way as in \cite[Lemma 2.7]{LZ2}.    \qed

\begin{lemma}\label{b3}
Let $(\mathcal{M}, \mathcal{Y})$ be a rigid cotorsion pair on $\C$. Then the following statements are equivalent.
\begin{itemize}
\item[{\rm (a)}] $\mathcal{M}$ is fully rigid.
\item[{\rm (b)}] $\underline{\C}\simeq \mathcal{H}/ \mathcal{M}$.
\item[{\rm (c)}] $\underline{\C}$ is abelian.
\end{itemize}
\end{lemma}

\proof ${\rm (a)}\Rightarrow {\rm (b)}$. It follows immediately from
\cite[Proposition 2.11]{LZ2}.

${\rm (b)}\Rightarrow {\rm (c)}$. By \cite[Theorem 3.2]{LN}, we have that $\mathcal{H}/ \mathcal{M}$ is abelian. Since  $\underline{\C}\simeq \mathcal{H}/ \mathcal{M}$, then $\underline{\C}$ is abelian.

${\rm (c)}\Rightarrow {\rm (a)}$. By Definition \ref{fullyrigid}, it is enough to show that any indecomposable object in $\C$ either belongs to $\mathcal{Y}$ or belongs to $\mathcal{H}$. Assume that $A$ is an indecomposable object and $A\notin \mathcal{Y}$. We will show $A\in \mathcal{H}$. Since $(\mathcal M,\mathcal Y)$ is cotorsion pair and (ET4)$^{\textrm{op}}$, we have a commutative diagram made of $\mathbb{E}$-triangles
$$\xymatrix{
Y_{1}\ar[r]\ar@{=}[d]&B\ar[r]^{m}\ar[d]^{h}&A\ar@{-->}[r]\ar[d]^{f}&\\
Y_{1}\ar[r]&M_{1}\ar[r]^{n}\ar[d]^{i}&Y\ar@{-->}[r]\ar[d]^{g}&\\
&M\ar@{=}[r]\ar@{-->}[d]&M\ar@{-->}[d]\\
&&}$$
where $M, M_1\in \mathcal{M}$ and $Y, Y_1\in \mathcal{Y}$. By definition, we have $B\in \mathcal{H}$. Next we show that $\underline{m}$ is an isomorphism. It suffices to show that $\underline{m}$ is both monic and epic, since $\underline{\C}$ is abelian.

First, we show that $\underline{m}$ is monic. If $\alpha:X\rightarrow B$ is a morphism in $\C$ such that $X$ has no direct summand in $\mathcal{Y}$ and $\underline{m\alpha}= 0$, then we obtain morphisms of $\mathbb{E}$-triangles as follows.
$$\xymatrix{X\ar[r]\ar[d]^{\alpha} &Y_2\ar[r] \ar[d]^k& M_2\ar@{-->}[r]\ar[d]^{j}&\\
B\ar[r]^{h} \ar[d]^{m} &M_1\ar[r]^{i} \ar[d]^n & M\ar@{-->}[r]\ar@{=}[d]&\\
A\ar[r]^{f} &Y \ar[r]^{g} & M \ar@{-->}[r]&}$$
where $M_2\in \mathcal{M}$ and $Y_2\in \mathcal{Y}$. The exactness of $\xymatrix{\C(Y_2, M_1)\ar[r]&\C(X, M_1)\ar[r]&\mathbb{E}(M_2, M_1)=0}$ shows the existence of $k$, the existence of $j$ is given by (ET3).  Since $m\alpha$ factors through $Y_2$, there is a morphism $l: M_2 \rightarrow Y$ such that $gl= j$ by \cite[Corollary 3.5]{NP}. Since $\mathbb{E}(M_2, Y_1)= 0$, there exists a morphism $p: M_2 \rightarrow M_1$ such that $l= np$. Hence $j= gl= gnp= ip$. By \cite[Corollary 3.5]{NP}, $\alpha$ factors through $Y_2$. Then $\underline{\alpha}= 0$.

Next, we show that $\underline{m}$ is epic. If $\beta:A\rightarrow Z$ is a morphism in $\C$ such that $\underline{\beta m}= 0$, then $\beta m$ factors through an object from $\mathcal{Y}$. Hence we obtain a commutative diagram
$$\xymatrix{
  B\ar[rr]^{\beta m} \ar[dr]_{y}
                &  &    Z     \\
                & Y\ar[ur]_{x}                 }$$
which means that $\beta m= xy$. Applying the functor $\C(-, Y)$ to the $\mathbb{E}$-triangle $$\xymatrix{B\ar[r]^{h}&M_{1}\ar[r]^{i}&M\ar@{-->}[r]&},$$ we obtain an exact sequence $$\C(M_{1}, Y)\xrightarrow{\C(h,Y)}\C(B, Y)\xrightarrow{~}\mathbb{E}(M, Y)=0.$$ Then for any morphism $y: B\rightarrow Y$, there exists a morphism $n\colon M_{1}\rightarrow Y$ such that $y=nh$. Thus $\beta m= xy= (xn)h$. By \cite[Lemma 3.13]{NP}, we have the following commutative diagram
$$
\xymatrix{
  B \ar[d]_{h} \ar[r]^{m}
                &A  \ar[d]_{f} \ar@/^/[ddr]^{\beta}  \\
  M_1 \ar[r]^{n} \ar@/_/[drr]_{xn}
                & Y \ar@{.>}[dr]|-{u}            \\
                &               & Z              }
$$
It follows that $\beta=uf$ which implies  $\underline{\beta}= 0$.

Since $A$ is indecomposable and $\underline{m}$ is an isomorphism, $A$ is a direct summand of $B$. By Lemma \ref{h0}, we have $A\in \mathcal{H}$.
This completes the proof.   \qed

\begin{lemma}\label{b4}
Let $\mathcal{M}$ be a fully rigid, $2$-rigid subcategory of $\C$. Then $\Omega (\mathcal{H}_{\mathcal{Y}}) \subseteq \mathcal{Y}$.
\end{lemma}

\proof Let $A\in \mathcal{H}_{\mathcal{Y}}$ be an indecomposable object admitting an $\mathbb{E}$-triangle $$\xymatrix{\Omega A\ar[r]^{g}&P_{A}\ar[r]^{f}&A\ar@{-->}[r]^{\delta}&}$$ where $P_{A}\in \mathcal{P}$ and $f$ is right minimal. If $\mathbb{E}(A, P_{A})\neq 0$, we have a non-split $\mathbb{E}$-triangle of the form
$\xymatrix{P_{A}\ar[r]^{q}&Q\ar[r]^{p}&A\ar@{-->}[r]^{\sigma}&}$. Since $A\in \mathcal{H}$, we have an $\mathbb{E}$-triangle $\xymatrix{A\ar[r]&M_{1}\ar[r]&M_{2}\ar@{-->}[r]^{\gamma}&}$ where $M_1, M_2\in \mathcal{M}$. Moreover, we can get the following commutative diagram
\begin{equation}\label{cd}
  \begin{array}{l}
$$\xymatrix{\Omega A\ar[r]^x\ar[d]^{g} &\Omega M_1\ar@{=}[r] \ar[d]^{n} &\Omega M_1 \ar[d]& &\\
P_A\ar[r]\ar[d]^{f} &\Omega M_2\ar[r]^{j} \ar[d]^m & P\ar[r]\ar[d]^{l}& M_2\ar@{-->}[r]\ar@{=}[d]&\\
A\ar@{=}[r]\ar@{-->}[d]^{\delta}&A\ar[r]^h\ar@{-->}[d]&M_1\ar@{=}[r]\ar@{-->}[d]&M_2\ar@{-->}[r]^{\gamma}&\\
&&&} $$
  \end{array}
\end{equation}
where $P, P_A\in \mathcal{P}$. Since $\mathcal{M}$ is 2-rigid, we have $\mathbb{E}(\Omega \mathcal{M}, \mathcal{M})= 0$. Then $\mathbb{E}(\Omega M_2, P_A)= 0$ since $P_A\in \mathcal{P}\varsubsetneq \mathcal{M}$. Hence, we get a morphism of $\mathbb{E}$-triangles.
$$\xymatrix{\Omega M_{1}\ar[r]^{n} \ar[d]^{b} &\Omega M_2\ar[r]^{m}\ar[d]^a & A\ar@{-->}[r]\ar@{=}[d]&\\
P_A\ar[r]^{q}&Q \ar[r]^{p} & A\ar@{-->}[r]&}$$
The exactness of $$\xymatrix{\C(\Omega M_2, Q)\ar[r]&\C(\Omega M_2, A)\ar[r]&\mathbb{E}(\Omega M_2, P_A)=0}$$ shows the existence of $a$, the existence of $b$ is given by (ET3)$^{\textrm{op}}$. Next we show that $\underline{p}$ is an epimorphism, it suffices to show that $\underline{m}$ is an epimorphism. If there is a morphism $e: A\rightarrow E$ such that $\underline{em}= 0$, then $em$ factors through an object from $\mathcal{M}$ since $\underline{\C}\simeq \mathcal{H}/ \mathcal{M}$. Hence it factors through $j$. Moreover, we have that $e$ factors through $M_1$, i.e. $\underline{e}= 0$. Then, we have $Q\notin \mathcal{Y}$. By \cite[Theorem 3.5]{LN}, $\underline{p}$ is monic in $\underline{\C}$, then $\underline{p}$ is an isomorphism. Since $A$ is indecomposable, $p$ is a retraction. Hence the $\mathbb{E}$-triangle $\xymatrix{P_{A}\ar[r]^{q}&Q\ar[r]^{p}&A\ar@{-->}[r]^{\sigma}&}$ is split, a contradiction. Thus, $\mathbb{E}(A, P_{A})= 0$. By Lemma \ref{b1}, $\Omega A$ is indecomposable. Since $\mathcal{M}$ is fully rigid, $\Omega A$ either belongs to $\mathcal{Y}$ or belongs to $\mathcal{H}$. If $\Omega A\in \mathcal{H}_\mathcal{Y}$, then there exists an $\mathbb{E}$-triangle $\xymatrix{\Omega A\ar[r]^{u}&M_1'\ar[r]^{v}&M_2'\ar@{-->}[r]^{\lambda}&}$ with $M_1'$ and $M_2'$ in $\mathcal{M}$. Then we obtain a morphism of $\mathbb{E}$-triangles
$$\xymatrix{\Omega A\ar[r]^{u} \ar@{=}[d] &M_1'\ar[r]^{v}\ar[d]^c & M_2'\ar@{-->}[r]^{\lambda}\ar[d]^d&\\
\Omega A\ar[r]^{g}&P_{A} \ar[r]^{f} & A\ar@{-->}[r]^{\delta}&}$$
where the existence of $c: M_1'\rightarrow P_{A}$ is given by the exactness of $$\xymatrix{\C(M_1', P_{A})\ar[r]&\C(\Omega A, P_{A})\ar[r]&\mathbb{E}(M_2', P_A)=0},$$ and the existence of $d: M_2'\rightarrow A$ is given by (ET3)$^{\textrm{op}}$. By the dual of \cite[Proposition 1.20]{LN}, we have an $\mathbb{E}$-triangle of the form $\xymatrix{M_1'\ar[r]&{M_2'} \oplus P_{A}\ar[r]&A\ar@{-->}[r]&}$. Applying the functor $\C(M_1', -)$ to the above $\mathbb{E}$-triangle, we obtain an exact sequence
$$\xymatrix{\mathbb{E}(M_1', {M_2'} \oplus P_{A})\ar[r]&\mathbb{E}(M_1', A)\ar[r]&\mathbb{E}^2(M_1', M_1')}.$$ Since $\mathcal{M}$ is 2-rigid, $\mathbb{E}(M_1', {M_2'}\oplus P_{A})= \mathbb{E}^2(M_1', M_1')=0$. Then $\mathbb{E}(M_1', A)= 0$, which implies $A\in \mathcal{M}^{\perp_1}= \mathcal{Y}$, this yields a contradiction. Hence $\Omega A \in \mathcal{Y}$. \qed
\medskip

This lemma immediately yields the following conclusion.

\begin{corollary}\label{b7}
Let $\mathcal{M}$ be a fully rigid, $2$-rigid subcategory of $\C$. Then $\underline{\C}$ is hereditary.

\proof It suffices to show that $\textrm{gl.dim} \underline{\C}\leq 1$. For any indecomposable object $A\in \mathcal{H}_{\mathcal{Y}}$, we have $\Omega A \in \mathcal{Y}$ by Lemma \ref{b4}. Then in the diagram (\ref{cd}), the morphism $x: \Omega A \rightarrow \Omega M_{1}$ factors through an object in $\mathcal{Y}$. By \cite[Proposition 2.12]{LZ2}, we have $\textrm{pd}_{\underline{\C} }(A)\leq  1$, as desired.    \qed
\end{corollary}

\begin{lemma}\label{b5}
Let $\mathcal{M}$ be a fully rigid subcategory of $\C$. For any indecomposable object $Y\in \mathcal{Y}_{\mathcal{M}}\cap {^{\perp_{1}}\mathcal{M}}$, one of the following conditions must be satisfied:
\begin{itemize}
\item[{\rm (a)}] $Y$ admits an $\mathbb{E}$-triangle $\xymatrix{\Omega Y\ar[r]&P\ar[r]^{f}&Y\ar@{-->}[r]&}$where $P \in \mathcal{P}$, $\Omega Y \in \mathcal{Y}$, and $f$ is right minimal.
\item[{\rm (b)}] $Y$ admits an $\mathbb{E}$-triangle $\xymatrix{M_2\ar[r]&M_1\ar[r]&Y\ar@{-->}[r]&}$where $M_1, M_2 \in \mathcal{M}$.
\end{itemize}
\end{lemma}

\proof For any indecomposable object $Y$, there exists an $\mathbb{E}$-triangle $$\xymatrix{\Omega Y\ar[r]&P\ar[r]^{f}&Y\ar@{-->}[r]&}$$ where $P \in \mathcal{P}$, $f$ is right minimal. Since $P \in \mathcal{P}\varsubsetneq \mathcal{M}$ and $Y\in \mathcal{Y}_{\mathcal{M}}\cap {^{\perp_{1}}\mathcal{M}}$, we have $\mathbb{E}(Y, P)=0$. By Lemma \ref{b1}, $\Omega Y$ is indecomposable. Since $\mathcal{M}$ is fully rigid, then $\Omega Y\in \mathcal{Y}$ or $\Omega Y\in \mathcal{H}$. In the first case $\Omega Y\in \mathcal{Y}$, then ${\rm (a)}$ holds. In the second case, by definition, we have an $\mathbb{E}$-triangle $\xymatrix{\Omega Y\ar[r]&M^1\ar[r]&M^2\ar@{-->}[r]&}$where $M^1, M^2 \in \mathcal{M}$. By \cite[Proposition 3.15]{NP}, we get the following commutative diagram made of $\mathbb{E}$-triangles
 $$\xymatrix{\Omega Y\ar[r]\ar[d]&M^1\ar[r] \ar[d]& M^2\ar@{-->}[r]\ar@{=}[d]&\\
P\ar[r]\ar[d] &M^3\ar[r] \ar[d] & M^2\ar@{-->}[r]^\theta&\\
Y\ar@{=}[r]\ar@{-->}[d]&Y\ar@{-->}[d]\\
&&}$$
Since  $\mathbb{E}(\mathcal{M}, \mathcal{M})=0$ and $\mathcal{P}\varsubsetneq \mathcal{M}$, we have $\theta=0$. It follows that
 $M^3\simeq M^2\oplus P\in\mathcal M$. This shows that the $\E$-triangle $\xymatrix{M^1\ar[r]& M^3\ar[r]&Y\ar@{-->}[r]&}$ is what we need. So ${\rm (b)}$ holds.  \qed
\medskip

From now on, let $\mathcal{M}$ be a $d$-cluster tilting subcategory in $\C$ with $d\geq 3$. In order to give our main theorem in this section, we need the following definition, which is the dual of \cite[Definition 5.4]{LN}.

\begin{definition}
For any $l\geq 0$, we define a subcategory $\mathcal{M}_l\subseteq \C$ inductively as follows.
\begin{itemize}
\item $\mathcal{M}_0= \mathcal{M}$.
\item $\mathcal{M}_l= \textrm{CoCone}(\mathcal{M}, \mathcal{M}_{l-1})$ for any $0< l\leq d-1$.
\end{itemize}
\end{definition}

\begin{remark}
From the above definition, we immediately see that $\mathcal{M}_1= \textrm{CoCone}(\mathcal{M}, \mathcal{M})= \mathcal{H}$, and $\mathcal{M}_0\subseteq \mathcal{M}_1 \subseteq\cdot \cdot \cdot \subseteq  \mathcal{M}_{d-1} =\C$. Moreover, by the dual of \cite[Proposition 5.10]{LN}, we have that $\mathcal{M}_l$ is closed under direct summands, for any $0\leq l\leq d-1$.
\end{remark}

For any $m>0$, we denote by $^{\perp_{m}}\mathcal{M}$ the subcategory of objects $X\in \C$ satisfying $$\mathbb{E}^{i}(X, \mathcal{M})=0\quad\mbox{for any}~~1\leq i\leq m.$$
We have $^{\perp_{d-1}}\mathcal{M}=\mathcal{M}$.

\begin{lemma}\label{zhengjiao}
By the dual of {\rm \cite[Proposition 5.2]{LN}}, we have the following result.
$$
\mathcal{M}_{l}= \left\{
\begin{array}{ll}
^{\perp_{d-l}}\mathcal{M}~~~ 1\leq l< d,\\[2mm]
\C~~~~~~~~~l\geq d.
\end{array}
\right.$$
\end{lemma}
\medskip

We are now in a position to prove our main result, which generalizes Liu's result \cite[Theorem 3.2]{L2} on an exact category.

\begin{theorem}\label{b6}
Let $\mathcal{M}$ be a $d$-cluster-tilting subcategory of $\C$ with $d\geq 3$. Assume that $(\mathcal{M}, \mathcal{Y}), (\mathcal{Y}, \mathcal{N})$ is the twin cotorsion pair induced by $\mathcal{M}$. Then
$$\mathcal{M} \mbox{ is fully rigid {\rm (}resp.}~ (\mathcal{M}, \mathcal{Y}), (\mathcal{Y}, \mathcal{N})\mbox{ has an abelian heart{\rm )} if and only if}~ ~\Omega(\mathcal{H}_{\mathcal{Y}})\subseteq \mathcal{Y}.$$ Moreover, this abelian heart $\underline{\C}$ is hereditary.
\end{theorem}

\proof {\bf Necessity.} It follows immediately from Lemma \ref{b4}.

 {\bf Sufficiency.} By Definition \ref{fullyrigid}, it suffices to show that any indecomposable object $A\notin \mathcal{Y}$ belongs to $\mathcal{H}_\mathcal{Y}$. We prove it by induction on $d$.

When $d= 3$, suppose that $A\in (\mathcal{M}_2)_\mathcal{Y}$ is any object, we have an $\mathbb{E}$-triangle $$\xymatrix{A\ar[r]&M^1\ar[r]&A^1\ar@{-->}[r]&}$$
where $M^1 \in \mathcal{M}$ and $A^1\in \mathcal{M}_1= \mathcal{H}$. If $A^1= A_1\oplus A_2$ such that the indecomposable object $A_1\in \mathcal{H}_\mathcal{Y}$ is non-zero, then we obtain morphisms of $\mathbb{E}$-triangles as follows.\\
  $$\xymatrix{\Omega A_1\ar[r]\ar[d]^{m} &P\ar[r]^f \ar[d]^a& A_1\ar@{-->}[r]\ar[d]^{\left[\begin{smallmatrix} 1\\[1mm] 0\end {smallmatrix}\right]}&\\
A\ar[r] \ar[d]^{n} &M^1\ar[r] \ar[d]^b & A_1\oplus A_2\ar@{-->}[r]\ar[d]^{\left[\begin{smallmatrix} 1& 0\end {smallmatrix}\right]}&\\
\Omega A_1\ar[r] &P \ar[r]^{f} & A_1 \ar@{-->}[r]&}$$
where $P\in \mathcal{P}$, $f$ is right minimal, and the existence of $b$ is given by the exactness of $$\xymatrix{\C(M^1, P)\ar[r]&\C(M^1, A_1)\ar[r]&\mathbb{E}(M^1, \Omega A_1)=0}.$$
By the definition of right minimal, we have $ba$ is an isomorphism.
By \cite[Corollary 3.6]{NP}, we obtain that $nm$ is an isomorphism. Since $A$ is indecomposable, we obtain $A\simeq \Omega A_1$. But $\Omega A_1\in \mathcal{Y}$, a contradiction. Thus $A^1\in \mathcal{M}$ and then $A\in \mathcal{H}_\mathcal{Y}$.

For $d> 3$, suppose that we have shown the case for $d= i$. Let $A\in (\mathcal{M}_{i})_\mathcal{Y}$ be any object, and take an $\mathbb{E}$-triangle $$\xymatrix{A\ar[r]^{g} &M^1\ar[r]&A^1\ar@{-->}[r]&}$$ where $A^1\in \mathcal{M}_{i-1}$. If $A^1= A_1'\oplus A_2'$ such that the indecomposable object $A_1'\notin \mathcal{Y}$, then by the assumption of the induction, $A_1'\in \mathcal{H}_\mathcal{Y}$. Since $\Omega (\mathcal{H}_\mathcal{Y})\subseteq \mathcal{Y}$, we obtain two morphisms of $\mathbb{E}$-triangles as follows.
$$\xymatrix{\Omega A_1'\ar[r]\ar[d]^{m^{\prime}} &P'\ar[r]^{f^{\prime}} \ar[d]^{a^{\prime}}& A_1'\ar@{-->}[r]\ar[d]^{\left[\begin{smallmatrix} 1\\[1mm] 0\end {smallmatrix}\right]}&\\
A\ar[r]^{g} \ar[d]^{n^{\prime}} &M^1\ar[r] \ar[d]^{b^{\prime}} & A_1'\oplus A_2'\ar@{-->}[r]\ar[d]^{\left[\begin{smallmatrix} 1& 0\end {smallmatrix}\right]}&\\
\Omega A_1'\ar[r] &P'\ar[r]^{f^{\prime}} & A_1' \ar@{-->}[r]&}$$
where $P'\in \mathcal{P}$ and $f'$ is right minimal. It is straightforward to verify that $A\simeq \Omega A_1'\in \mathcal{Y}$, this yields a contradiction.
Thus $A^1\in \mathcal{Y}$.
Since $A^1\in\mathcal M_{i-1}$, we have $A^1\in{^{\perp_{i-(i-1)}}\mathcal{M}}={^{\perp_{1}}\mathcal{M}}$
 by Lemma \ref{zhengjiao}.
It follows that $A^1\in \mathcal{Y}_{\mathcal{M}}\cap{^{\perp_{1}}\mathcal{M}}$. If $A^1= A_1''\oplus A_2''$ such that the indecomposable object $A_1''$ satisfying the condition $(a)$ in Lemma \ref{b5}, we obtain two morphisms of $\mathbb{E}$-triangles as follows.
$$\xymatrix{\Omega A_1''\ar[r]\ar[d]^{m^{\prime \prime}} &P''\ar[r]^{f^{\prime \prime}} \ar[d]^{a^{\prime \prime}}& A_1''\ar@{-->}[r]\ar[d]^{\left[\begin{smallmatrix} 1\\ 0\end {smallmatrix}\right]}&\\
A\ar[r] \ar[d]^{n^{\prime \prime}} &M^1\ar[r] \ar[d]^{b^{\prime \prime}} & A_1''\oplus A_2''\ar@{-->}[r]\ar[d]^{\left[\begin{smallmatrix} 1& 0\end {smallmatrix}\right]}&\\
\Omega A_1''\ar[r] &P''\ar[r]^{f^{\prime \prime}} & A_1'' \ar@{-->}[r]&}$$
where $P''\in \mathcal{P}$ and $f''$ is right minimal. We can obtain $A\simeq \Omega A_1'' \in \mathcal{Y}$, a contradiction. By Lemma \ref{b5}, $A^1$ admits an $\mathbb{E}$-triangle $\xymatrix{M_2\ar[r]&M_1\ar[r]&A^1\ar@{-->}[r]&}$where $M_1, M_2 \in \mathcal{M}$.
By \cite[Proposition 3.15]{NP}, we get the following commutative diagram made of $\mathbb{E}$-triangles
  $$\xymatrix{&A\ar@{=}[r] \ar[d]&A \ar[d]\\
M_2\ar[r]  \ar@{=}[d] &M\ar[r] \ar[d] & M^1\ar[d]\ar@{-->}[r]^{\;\; \omega}&\\
M_2\ar[r]&M_1\ar[r]\ar@{-->}[d]&A^1\ar[r]\ar@{-->}[d]&\\
&&}$$
Since $M^1,M_2\in\mathcal M$ and $\mathbb{E}(\mathcal{M}, \mathcal{M})=0$, we have $\omega=0$ which implies $M\simeq M_2\oplus M^1\in\mathcal M$.
Thus the $\E$-triangle
$\xymatrix{A\ar[r]&M\ar[r]&M_1\ar@{-->}[r]&}$ shows that $A\in \mathcal{H}_\mathcal{Y}$.
\smallskip

It remains to show $\underline{\C}$ is hereditary. It follows immediately from Corollary \ref{b7}.     \qed

\begin{remark}
In Theorem \ref{b6}, when $\C$ is an exact category, it is just the Theorem 3.2 in \cite{L2}, and
when $\C$ is a triangulated category, it is a new phenomena.
\end{remark}

\textbf{Qiong Huang}\\
School of Mathematics and Statistics,~ Hunan Normal University,
410081 Changsha, Hunan, P. R. China.\\
E-mail: qhuang@hunnu.edu.cn\\[2mm]
\textbf{Panyue Zhou}\\
College of Mathematics, Hunan Institute of Science and Technology, 414006 Yueyang, Hunan, P. R. China.\\
E-mail: panyuezhou@163.com

\end{document}